\documentclass[12pt]{article}
\usepackage{amsmath}
\begin{document}
\title{An easier solution of a Diophantine problem about triangles,
in which those lines from the vertices which bisect the opposite sides may
be expressed rationally\footnote{Delivered to the St.-Petersburg
Academy on August 12, 1779. Originally published as
\emph{Solutio facilior problematis Diophantei circa triangulum, in quo rectae
ex angulis latera opposita bisecantes rationaliter exprimantur},
M\'emoires de l'Acad\'emie Imp\'eriale des Sciences de St.-P\'etersbourg
\textbf{2} (1810),
10-16, and reprinted in \emph{Leonhard Euler, Opera Omnia}, Series 1:
Opera mathematica, Volume 4, Birkh\"auser, 1992.
A copy of the original text is available
electronically at the Euler Archive, at www.eulerarchive.org. This paper
is E732 in the Enestr\"om index.}}
\author{Leonhard Euler\footnote{Date of translation: March 2, 2005.
Translated from the Latin by Jordan
Bell, 2nd year undergraduate in Honours Mathematics, School of Mathematics
and Statistics, Carleton University, Ottawa, Ontario, Canada. Email:
jbell3@connect.carleton.ca. Translation written under Dr. B. Mortimer.}}
\date{}

\maketitle

1. The sides of a triangle having been placed as 
$AB=2c$, $AC=2b$ and $BC=2a$, if the bisecting lines are called
$AX=x$, $BY=y$ and $CZ=z$, it is apparent from Geometry 
for the squares of these three lines to be expressed in the following way: 
\begin{align*}
xx=2bb+2cc-aa\\
yy=2cc+2aa-bb\\
zz=2aa+2bb-cc,
\end{align*}
and these equations are thus to be resolved through rational numbers.

2. From these three equations, the following three are formed:
\begin{align*}
\textrm{I.} \quad &xx-yy=3(bb-aa)\\
\textrm{II.} \quad &xx+yy=4cc+aa+bb\\
\textrm{III.} \quad &zz=2aa+2bb-cc.
\end{align*}
We deal with these equations in the following way.

3. We start with the first of these equations, and so that we can
avoid fractions, we set $xx-yy=3(bb-aa)=12fg(pp-qq)$; from this,
with 
it  $bb-aa=4fg(pp-qq)$, we make it $b+a=2f(p+q)$ and
$b-a=2g(p-q)$, from which it will be $b=(f+g)p+(f-g)q$ and
$a=(f-g)p+(f+g)q$. Indeed, with the previous equation
$xx-yy=12fg(pp-qq)$, we take
$x+y=6g(p+q)$ and $x-y=2f(p-q)$, from which it will be
$x=(3g+f)p+(3g-f)q$ and $y=(3g-f)p+(3g+f)q$.

4. The second equation is now approached:
\[ xx+yy=4cc+aa+bb, \]
and together with the values that have just been determined, it is
 discovered that:
\[ xx+yy=2(9gg+ff)(pp+qq)+4(9gg-ff)pq. \]
 Then indeed it will be:
\[ bb+aa=2(ff+gg)(pp+qq)+4(ff-gg)pq.\]
 Therefore, with it $4cc=xx+yy-(bb+aa)$, it will be:
\[ 4cc=16gg(pp+qq)+(40gg-8ff)pq, \]
on account of which we have $cc=4gg(pp+qq)+(10gg-2ff)pq$,
so that therefore this formula is to be transformed into a square.

5. From here the third equation is proceded to, which is
$zz=2(aa+bb)-cc$, from which it will be:
\begin{align*}
2(aa+bb)=4(ff+gg)(pp+qq)+8(ff-gg)pq \quad \textrm{ and}\\
cc=4gg(pp+qq)+(10gg-2ff)pq, \quad \textrm{ and it will be}\\
zz=4ff(pp+qq)+(10ff-18gg)pq,
\end{align*}
which is thus another formula that we must reduce to a square.

6. We divide the former of these two equations by $4gg$
and the latter indeed by $4ff$, so that we may have the following
formulas, reduced into squares:
\begin{align*}
\frac{cc}{4gg}=pp+qq+2pq \Big ( \frac{5gg-ff}{4gg} \Big ) \\
\frac{zz}{4ff}=pp+qq+2pg \Big ( \frac{5ff-9gg}{4ff} \Big ).
\end{align*}

7. For briefness, here we put $\frac{5gg-ff}{4gg}=m$
and $\frac{5ff-9gg}{4ff}=n$, so that now the entire problem
is brought to the resolution of these two formulas:
\begin{align*}
\frac{cc}{4gg}=pp+qq+2mpq=tt\\
\frac{zz}{4ff}=pp+qq+2npq=uu,
\end{align*}
and with this we will have $c=2gt$ and $z=2fu$. At this point,
in order for these two formulas to be reduced to squares, it is noted to
be that $tt-uu=2(m-n)pq$, and for convenience it is thus permitted
to treat this such that it is set $t+u=(m-n)p$ and $t-u=2q$, from which
is gathered $t=\frac{1}{2}(m-n)p+q$ and $u=\frac{1}{2}(m-n)p-q$.

8. But if these values are now substituted in place of $t$ and $u$,
 it comes about that squares simultaneously arise on both sides of the equation:
\[ pp\Big ( 1-\frac{1}{4}(m-n)^2 \Big )+(m+n)pq=0,\]
and when this equation is divided by $p$ it is:
\[ p\Big ( 1- \frac{1}{4}(m-n)^2 \Big ) + (m+n)q=0,\]
from which the ratio between the 
letters $p$ and $q$ freely comes forth: $\frac{p}{q}=\frac{4(m+n)}{(m-n)^2-4}$,
on account of which it will be possible for it to be taken
$p=4(m+n)$ and $q=(m-n)^2-4$, or equivalent multiples, so that
 it is to be held in general that $p=4(m+n)M$ and
$q=((m-n)^2-4)N$, and in this way all the conditions are fully
satisfied.

9. Now with the letters $p$ and $q$ that have been worked out, it will
be $t=2(m+n)(m-n)+(m-n)^2-4=(m-n)(3m+n)-4$ and $u=(m-n)(m+3n)+4$,
from which as before we infer these determinations:
\begin{align*}
c=2g(m-n)(3m+n)-8g \quad \textrm{ and}\\
z=2f(m-n)(m+3n)+8f.
\end{align*}
Moreover, it is helpful to sum in the $p$ and $q$ places
the values of them from \S 7, so that it will thus be 
$c=g(m-n)p+2gq$ and $z=f(m-n)p-2fq$.

10. Now therefore a solution of our problem is able to be constructed
in the following way:

1.) Both the letters $f$ and $g$ are allowed to be at our choice, from
which
the letters $m$ and $n$ may be defined by means of the formulas
$m=\frac{5gg-ff}{4ff}$ and $n=\frac{5ff-9gg}{4ff}$.

2.) Now as before the letters $p$ and $q$ are sought through the power
of the formulas $p=4(m+n)$ and $q=(m-n)^2-4$, with which having been
obtained, the lines bisecting the sides of the triangle can be
expressed in the following way.

3.) Of course the sides are found with these formulas:
\begin{align*}
a=(f-g)p+(f+g)q\\
b=(f+g)p+(f-g)q\\
c=2g(m-n)(3m+n)-8g=g(m-n)p+2gq
\end{align*}

4.) Finally therefore, the bisecting lines themselves will be had:
\begin{align*}
x=(3g+f)p+(3g-f)q\\
y=(3g-f)p+(3g+f)q\\
z=2f(m-n)(m+3n)+8f=f(m-n)p-2fg.
\end{align*}

11. So that we can illustrate an example of this,
we suppose that $f=2$ and $g=1$, and it will then be $m=\frac{1}{4}$
and $n=\frac{11}{16}$. From this as before it is gathered that
$p=\frac{15}{4}$ and $q=-\frac{975}{256}$, whose values, to be
reduced to integral numbers, are given as $p=64$ and $q=-65$. Then
our triangle will now be determined in this way:
\begin{align*}
a=131; b=127; c=158; \quad \textrm{ then indeed,}\\
x=255; y=261; z=204.
\end{align*}

12. On this occasion as has been noted for the example,
it will be helpful now too for the letters $x,y,z$ to be taken for
the sides, with in general it being:
\begin{align*}
2xx-2yy-zz=9cc\\
2yy+2zz- -xx=9aa\\
2zz+2xx-yy=9bb,
\end{align*}
from which it follows, with it permitted for
the numbers yielding $x,y,z$ to be reduced by 3, then a very simple
triangle is able to be formed, of which the sides are 136, 170, 174.

13. At this point, with it $m=\frac{5gg-ff}{4gg}$ and $n=\frac{5ff-9gg}{4ff}$,
it will first be $m+n=\frac{10ffgg-f^4-9g^4}{4ffgg}=
-\frac{(gg-ff)(9gg-ff)}{4ffgg}$. Then indeed it will be had
that $m-n=\frac{9g^4-f^4}{4ffgg}=\frac{(3gg+ff)(3gg-ff)}{4ffgg}$,
from which it is gathered that $m-n+2=\frac{9g^4+8ffgg-f^4}{4ffgg}=
\frac{(gg+ff)(9gg-ff)}{4ffgg}$ and
$m-n-2=\frac{9g^4-8ffgg-f^4}{4ffgg}=\frac{(gg-ff)(9gg+ff)}{4ffgg}$.
By means of this, with it having been found that $p=4(m+n)M$, it will
now be:
\begin{align*}
p=-\frac{(gg-ff)(9gg-ff)M}{4ffgg} \quad \textrm{ and,}\\
q=((m-n)^2-4)M=\frac{(g^4-f^4)(81g^4-f^4)}{16f^4g^4}M.
\end{align*}

14. Therefore we now reduce the ratio of the letters $p$ and $q$ to minimum
terms, and we take $M=\frac{16f^4g^4}{(gg-ff)(9gg-ff)}$,
and thus it comes forth that $p=-16ffgg$ and $q=(gg+ff)(9gg+ff)$. From
this, because it was $t=\frac{1}{2}(m-n)p+q$, it now follows that,
\[ t=(gg+ff)(9gg+ff)-2(3gg+ff)(3gg-ff). \]
In a similar way it will be $u=\frac{1}{2}(m-n)p-q$, from which it is gathered
that:
\[ u=-(gg+ff)(9gg+ff)-2(3gg+ff)(3gg-ff). \]
Then because of this it will be $c=2gt$ and $z=2fu$.

15. For the remaining letters, it will first be $a+b=2f(p+q)$ and
$b-a=2g(p-q)$; then indeed it will be $x+y=6g(p+q)$ and $x-y=2f(p-q)$.
It is no problem for however many examples to be resolved
readily enough by means of these formulas, to such a degree as they
are according to integers; since they are in every respect determined
by the ratio of the integral numbers $f$ and $g$, these are therefore
always allowed to be assumed integers.

16. It will be $f=1$ and $g=2$, and it will be $p=-64$ and $q=185$.
Hence now it is gathered $t=-101$ and $u=-471$; then indeed $c=-404$
and $z=-942$. To be sure, as before it will be $b-a=-996$ and
$b+a=+242$; $x+y=1452$ and $x-y=-498$. It will therefore be:
\begin{align*}
a=619; b=377; c=404\\
x=477; y=975; z=942.
\end{align*}
It is observed that the numbers $x,y,z$ are able moreover
to be used instead of
$a,b,c$, with these divided by 3, so that it becomes:
\begin{align*}
a=159; b=325; c=314, \quad \textrm{ then also it will be,}\\
x=619; y=377; z=404.
\end{align*} 

\end{document}